\documentclass[a4paper]{article}

\usepackage[english]{babel}
\usepackage[utf8x]{inputenc}
\usepackage[T1]{fontenc}

\usepackage[a4paper,top=3cm,bottom=2cm,left=3cm,right=3cm,marginparwidth=1.75cm]{geometry}

\usepackage{amsmath}
\usepackage{graphicx}
\usepackage{bm}
\usepackage{caption}
\usepackage{float}
\usepackage{amssymb}
\usepackage{graphicx}
\usepackage[colorinlistoftodos]{todonotes}
\usepackage[colorlinks=true, allcolors=blue]{hyperref}
\usepackage[numbers,sort&compress]{natbib}
\usepackage{amsthm}
\usepackage[title]{appendix}
\usepackage{indentfirst} 
\usepackage{authblk}
\usepackage{makecell}
\usepackage{bm}
\usepackage{amsmath}
\usepackage{chngcntr}

\numberwithin{figure}{section}

\theoremstyle{plain}

\theoremstyle{definition}

\theoremstyle{remark}

\counterwithout{figure}{section} 
\counterwithout{table}{section} 

\numberwithin{equation}{section}

\begin{document}



\title{\bf\Large A High-resolution Spatiotemporal Coupling Ghost Fluid Method for Two-Dimensional Compressible Multimedium Flows with Source Terms}
\author[1,2]{Zhixin Huo\thanks{Corresponding author. Email addresses:
zhixinhuo@hpu.edu.cn (Z. Huo)
}}

\affil[1]{School of Mathematics and Information Science, Henan Polytechnic University, Jiazuo, Henan, 454000, PR China}
\affil[2]{School of Mechatronical Engineering, Beijing Institute of Technology, Beijing, 100081, PR China}
\date{}

\maketitle

\begin{abstract}


While exact and approximate Riemann solvers are widely used, they exhibit two fundamental limitations:
1) 
Fail to represent continuous entropy transport processes, resulting in thermodynamic incompatibility that limits their applicability to compressible flows.
2) Consider only the effects of normal components at interfaces while neglecting the effects of tangential flux and source term, making them unsuitable for multidimensional problems and cases involving source terms. These limitations persist in Riemann problem-based ghost fluid methods.
To address these challenges, we developed a novel spatiotemporal coupling high-resolution ghost fluid method featuring two key advancements:
1) Integration of nonlinear geometrical optics to properly account for thermodynamic entropy evolution.
2) Implementation of the Lax-Wendroff/Cauchy-Kowalevski approach to incorporate tangential fluxes and source term effects.
These enhancements have been systematically applied to Riemann problem-based ghost fluid methods. Comprehensive numerical experiments demonstrate significant improvements in simulation accuracy and robustness compared to conventional approaches.


\end{abstract}

\section{Introduction}

The ghost fluid method (GFM) fundamentally operates by decomposing multi-medium problems into single-medium subproblems through the strategic definition of ghost fluid  regions and corresponding ghost fluid states. The distinguishing feature among various GFM implementations lies in their respective approaches to defining these ghost fluid states.
The seminal GFM formulation \cite{Fedkiw-1, Fedkiw-2} employed interpolation techniques for ghost fluid states determination, though this approach demonstrated limited applicability. Subsequent developments addressed specific challenges: Liu et al.\cite{Liu-2003}  introduced a modified GFM to handle strong interfacial waves, followed by a specialized water-gas interface treatment \cite{ref9}. Wang et al.\cite{Wang-2006} subsequently developed the real GFM to resolve vacuum-related complications. These advances precipitated numerous derivative methods, including the interface interaction GFM \cite{Hu-2004}, practical GFM \cite{Xu-2016}, and the MGFM/AC approach  \cite{Liu-Feng-2019}.
While these Riemann-problem-based GFMs incorporate material properties and interfacial interactions, they exhibit two critical limitations that can be addressed through generalized Riemann problem (GRP) theory:
1) Failure to properly account for thermdynamic entropy evolution restricts applicability to compressible flows;
2) Neglect the effects of tangential flux and source term precludes accurate multidimensional and source-term-included simulations.
The remainder of this paper first establishes the theoretical framework of GRP, then systematically addresses these two fundamental issues through detailed analysis and methodological improvements.

The differences between generalized Riemann problem and Riemann problem lie in the following two aspects: 
1) The initial values of the Riemann problem are constant on both sides of the discontinuity, while those of the generalized Riemann problem are not constant on both sides of the discontinuity.
2) The solution of the Riemann problem does not take into account the source term, while the solution of the generalized Riemann problem does take the source term into consideration.
Menshov \cite{Men'shov-1990}, Kolgan\cite{Kolgan-1972}, van Leer\cite{Leer-2006} et al. studied this kind of problems.
The Generalised Riemann Problem, or GRP, method was first  applied to compressible flow
by Ben-Artzi and Falcovitz in  \cite{Ben-Artzi-1984, Ben-Artzi-1985, Ben-Artzi-1986,Ben-Artzi-2003}. Ben-Artzi and Jiequan Li  proposed the GRP solver directly in Euler scheme and apply it to hyperbolic conservation laws in \cite{Ben-Artzi-2006, Ben-Artzi-2007}.
Then extended the GRP solver to high order in \cite{Qian-2014,Wu-2014, Wang-2015}.

Whether there are thermodynamic effects is the key to distinguishing whether a fluid is compressible or not. Moreover, the stronger the compressibility of a fluid is, the more important the thermodynamic effects becomes. The Gibbs relation is the core tool for thermodynamic analysis, and its mathematical expression is:
\begin{equation}
Tds=de-\frac{p}{\rho^2}d\rho,
\label{Gibbs}
\end{equation}
where $T$ is the temperature, $s$ is the entropy, $e$ is the internal energy, $p$ is the 
pressure and $\rho$ is the density.
The initial data of Riemann problem is piecewise constant, which implies that the
solutions emanating from the jumps are self-similarity and in particular rarefaction waves are always isentropic. Although the shock wave caused by the jump is non-isentropic, the Riemann problem can only depict the sudden entropy change but cannot reflect the continuous entropy transport process.  Therefore, according to the Gibbs relation \eqref{Gibbs},
the Riemann problem is thermodynamically incompatible, and is not applicable for compressible flows. While the initial data of generalized Riemann problem (GRP) is piecewise
polynomials of high degree, which implies that all waves from the initial jumps are curved and in particular rarefaction waves are no longer isentropic. Then the interaction between 
the kinematical and thermodymical quantities are activated by the initial entropy variation. By making use of the tricky technique of "nonlinear geometric optics",
that is, by introducing the transformation between spatial-temporal coordinates  and 
characteristic coordinates, to the singularity point, 
the spatial-temporal couping GRP solver can plug the entropy variation into numerical fluxes so as to precisely characterize the thermodynamical process, as described in \cite{GRP-2017}.   

The numerical fluxes which based on exact or approximate Riemann solvers only consider 
the Riemann problam along the normal direction, but neglect the the effects of tangential fluxes and source terms, which is not compatible to multi-dimensional problems and those with source terms. By making use of the Lax-Wendroff or Cauchy-Kowalevski method, that is, the unknown time derivatives of physical quantities can be expressed by the known spatial derivatives of physical quantities, the spatial-temporal GRP solver can plug the effects of tangential fluxes \cite{XinLei-2019} and source terms \cite{ref20, Zhu-2023, Huo-2022} into the numerical fluxes.

The Riemann problem based GFM define constantly distributed ghost fluid states in ghost fluid regions,
while the GRP based GFM our proposed define linely distributed ghost fluid states in ghost fluid regions. Here the spatial derivatives of the ghost fluid states along the normal direction of the interface is determined based on the post-wave states of the local double-medium generalized Riemann problem established along the normal direction of the interface. The ghost fluid states of each medium is determined by the corresponding post-wave states of such medium. The spatial derivatives of the ghost fluid states along the tangential direction is unchanged. 

In this paper, we consider the compressible multi-medium problems with source term in two dimensions,  in order to reflect the thermodynamic effects of compressible fluids, as well as the effects of tangential and source terms, we use the  technique of "nonlinear geometric optics" 
and the  Cauchy-Kowalevski method, and replace the numerical flux and ghost fluid methods based on the Riemann problem with those based on the generalized Riemann problem. In the second section we 
present the mathematical model of two-dimensional compressible multimedium flow problems. 
Introduce the general framework of the ghost fluid method used for handling multi-medium problems and the general framework of the finite volume scheme used for dealing with single-medium problems.
In the third section,  we introduce the Riemann problem based numerical fluxes and ghost fluid method.
In the fourth section, we introduce the Riemann problem based numerical fluxes and ghost fluid method.
In the fifth section, for the two-dimensional compressible flow problem with axisymmetric source terms, we have provided several typical numerical examples. The numerical results have demonstrated the superiority of the new algorithm proposed by us.

\section{The two-dimensional multi-medium compressible flow problems with source term}

\begin{figure}[H]
  \centering
  \includegraphics[width=8cm]{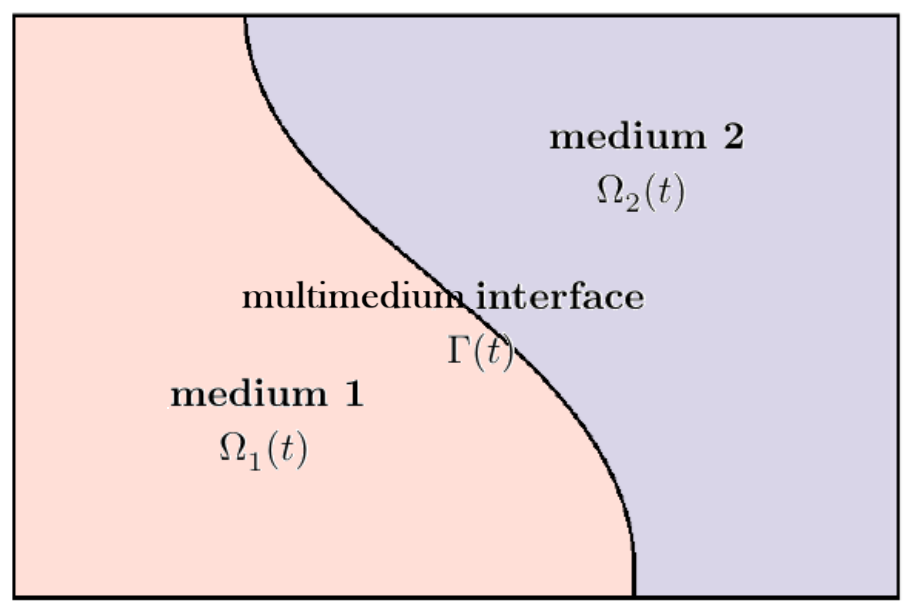}\\
  \caption{The mathematical model of the two-dimensional compressible multi-medium flow.}
  \label{MathModel}
\end{figure}

The two-dimensional compressible multi-medium flow problem with source terms is as shown in Figure \ref{MathModel}, the mathematical model of which is as follows:

\begin{equation}
\begin{aligned}
&\frac{\partial\bm{u}}{\partial t}+\frac{\partial\bm{F}}{\partial x}(\bm{u})+\frac{\partial\bm{G}}{\partial y}(\bm{u})={\bm {S}}(x,y,\bm{u}), \ \ (x,y)\in\Omega, \ t>0,\\
&{\bm{u}}(x,y,0)=
\begin{cases}
{\bm{u}}_1(x,y),\ \ (x,y)\in \Omega_1^0,\\
{\bm{u}}_2(x,y),\ \ (x,y)\in \Omega_2^0,
\end{cases}
\quad e(\rho,p)=
\begin{cases}
e_1(\rho,p),\ \ (x,y)\in \Omega_1(t),\\
e_2(\rho,p),\ \ (x,y)\in \Omega_2(t),
\end{cases}
\end{aligned}
\end{equation}
where
$$
{\bm{u}}=
\left[
\begin{matrix}
\rho\\
\rho u_x\\
\rho u_y\\
E
\end{matrix}
\right],\quad
{\bm{F}}=\left[
\begin{matrix}
\rho u_x\\
\rho u_x^2+p\\
\rho u_x u_y\\
u(E+p)
\end{matrix}
\right],\quad
{\bm{G}}=
\left[
\begin{matrix}
\rho u_y\\
\rho u_x u_y\\
\rho u_y^2+p\\
v(E+p)
\end{matrix}
\right],
$$
$\rho$ is the density,
$V=(u_x,u_y)$ is the velocity,
$p$ is the pressure,
$e$ is the specific internal energy,
$E=e+\frac{1}{2}(u_x^2+u_y^2)$ is the specific energy,
$\Omega$ is the computation domain, 
$\Gamma(t)$ is the multimedium interface at the time $t$,
$\Omega_k(t)$ is the domain of medium $k$ at the time $t$,
 $\Omega_k^0\triangleq\Omega_k(t_0)$, 
$e=e_{k}(\rho,p)$ is the EOS of medium $k$,
$k=1,2,$
$\Omega=\Omega_1(t)\cup\Gamma(t)\cup\Omega_2(t)$,
and ${\bm{S}}(x,y,{\bm{u}})$ is the source term.

\subsection{The ghost fluid method for multi-medium flow problems}

\begin{figure}[H]
  \centering
  \includegraphics[width=10cm]{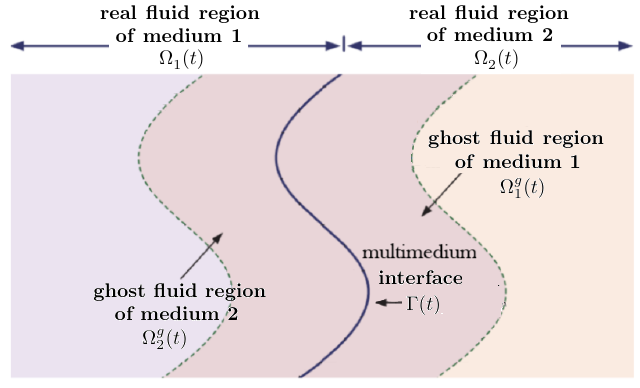}\\
  \caption{The schedule of the ghost fluid method for double medium in 2D.}
  \label{ghost_fluid_method_in_2D}
\end{figure}

The strategy of ghost fluid method (GFM) is that by introducing ghost fluid regions and ghost fluid states,
the multimedium problem can be decoupled into several single medium problems to solve separately.
As shown in Figure \ref{ghost_fluid_method_in_2D},
we define $\Omega_{k}^{g} (t)$ as the ghost fluid regions of medium $k$ at the time $t$, $k=1,2$, which are
Narrowband regions near the multimedium interface  in the regions of another medium, i.e., $\Omega_{1}^{g}(t)\subset\Omega_{2}(t)$, $\Omega_{2}^{g}(t)\subset\Omega_{1}(t)$. Then we define the ghost fluid states on the ghost fluid regions at the time $t=t_{n}$, i.e.,
\begin{equation}
\bm{u}(x,y,t_n)=\bm{u}_{k}^{g,n}(x,y), \quad (x,y)\in\Omega_{k}^{g,n},\; k=1,2,
\end{equation}
where $\Omega_{k}^{g,n}\triangleq\Omega_{k}^{g}(t_{n})$. 
The differences between different ghost fluid methods are
the definition of the ghost fluid states. The Riemann problem based ghost fluid method
defines constantly distributed ghost fluid states, while
The GRP based ghost fluid method defines linearly distributed ghost fluid states. 
Compared to the former, the latter has higher precision and can also reflect the thermal hydrodynamic effects of compressible fluids, as well as the effects of tangential flux and source terms, which will be  elaborated in the text later.

We solve the single medium problem in the real fluid region and ghost fluid region for each medium, i.e.,
\begin{equation}
\begin{aligned}
&\frac{\partial\bm{u}_{k}}{\partial t}+\frac{\partial \bm{F}}{\partial x}(\bm{u}_k)+\frac{\partial\bm{G}}{\partial y} (\bm{u}_k)={\bm {S}}(x,y,\bm{u}_k), \ \ (x,y)\in\Omega_k^n\cup\Gamma^n\cup\Omega_k^{g,n}, \ t>0,\\
&{\bm{u}}(x,y,t_n)=
\begin{cases}
{\bm{u}}_k^{n}(x,y),& (x,y)\in \Omega_k^n,\\
{\bm{u}}_{k}^{g,n}(x,y),& (x,y)\in \Omega_k^{g,n},
\end{cases}
\quad e(\rho,p)=e_k(\rho,p),\ \ (x,y)\in \Omega_k^n\cup\Gamma^n\cup\Omega_{k}^{g,n},\\
&k=1,2.
\end{aligned}
\label{SingleMediumProblem}
\end{equation}
We solve the single medium problem \eqref{SingleMediumProblem} for 
\begin{equation}
\bm{u}_k(x,y,t_{n+1})=\bm{u}_k^{n+1}(x,y),\ \ (x,y)\in \Omega_k^n\cup\Gamma^n\cup\Omega_{k}^{g,n},\ k=1,2.
\label{UpdataSingleMediumFluidStates}
\end{equation}

We update the multimedium interface $\Gamma(t)$ by the Level-Set method, i.e.,
solve the problem 
\begin{equation}
\begin{aligned}
&\frac{\partial\phi}{\partial t}+u_x\frac{\partial\phi}{\partial x}+u_y\frac{\partial \phi}{\partial y}=0, \ \ (x,y)\in\Omega,\ t>0,\\
&\phi(x,y,t_n)=\phi^{n}(x,y),\ \ (x,y)\in\Omega,
\end{aligned}
\end{equation}
where $\phi$ is the Level-Set function and $(u,v)$ is the vector of velocity, then the new location of the multimedium interface
$\Gamma^{n+1}$ as well as the new location of the real fluid region of medium $k$, i.e., $\Omega_k^{n+1}$, $k=1,2$, are defined by
\begin{equation}
\begin{aligned}
&\Omega_{1}^{n+1}:=\{(x,y)\in\Omega\lvert\phi^{n+1}(x,y)<0\},\\
&\Gamma^{n+1}:=\{(x,y)\in\Omega\lvert\phi^{n+1}(x,y)=0\},\\
&\Omega_{2}^{n+1}:=\{(x,y)\in\Omega\lvert\phi^{n+1}(x,y)>0\}.
\end{aligned}
\label{UpdateMultimediumInterface}
\end{equation}

Then we update the real fluid states in the whole computation region by combining \eqref{UpdataSingleMediumFluidStates} and \eqref{UpdateMultimediumInterface}, i.e.,
\begin{equation}
{\bm{u}(x,y,t_{n+1})}=
\begin{cases}
{\bm{u}}_1^{n+1}(x,y),\ \ (x,y)\in\Omega_1^{n+1},\\
{\bm{u}}_2^{n+1}(x,y),\ \ (x,y)\in\Omega_2^{n+1}.
\end{cases}
\end{equation}
The time step is updated by CFL condition, i.e., 
\begin{equation}
\Delta t_n=C_{\rm CFL}\frac{\min\{\Delta x,\Delta y\}}{\max_{i,j}\sqrt{{(u_{i,j}^n)^2+(v_{i,j}^n)^2}}}.
\end{equation}
We suggest the computation domain is $\Omega=[a,b]\times[c,d]$, and the teminate computation time is $T$, and we divide the space and time as follows:
$x_{\frac{1}{2}}=a$, $x_{i+\frac{1}{2}}=x_{i-\frac{1}{2}}+\Delta x$, $i=1,2,\cdot\cdot\cdot,\frac{b-a}{\Delta x}$;
$y_{\frac{1}{2}}=c$, $y_{i+\frac{1}{2}}=y_{i-\frac{1}{2}}+\Delta y$, $j=1,2,\cdot\cdot\cdot,\frac{d-c}{\Delta y}$;
$t_{0}=0$, $t_{N}=T$, $t_{n+1}=t_{n}+\Delta t_n$, $n=0,1,\cdot\cdot\cdot, N-1$.

\subsection{The finite volume method for single-medium flow problems}

Integrate the single-medium flow problems \eqref{SingleMediumProblem} on the region $[t_n,t_{n+1}]\times[x_{i-\frac{1}{2}},x_{i+\frac{1}{2}}]\times[y_{j-\frac{1}{2}},y_{j+\frac{1}{2}}]$, 
and we get the finite volume method for \eqref{SingleMediumProblem} as follows,
\begin{equation}
\begin{aligned}
&\frac{\overline{\bm{u}}_{i,j}^{n+1}-\overline{\bm{u}}_{i,j}^n}{\Delta t}\\
&+\frac{\frac{1}{\Delta t}\int_{t_n}^{t_{n+1}}{\bm{F}(\bm{u}(x_{i+\frac{1}{2}},y_j,t))}dt-\frac{1}{\Delta t}\int_{t_n}^{t_{n+1}}{\bm{F}(
\bm{u}(x_{i-\frac{1}{2}},y_j,t))}dt}{\Delta x}\\
&+\frac{\frac{1}{\Delta t}\int_{t_n}^{t_{n+1}}{\bm{G}(\bm{u}(x_i,y_{j+\frac{1}{2}},t))}dt-\frac{1}{\Delta t}\int_{t_n}^{t_{n+1}}{\bm{G}(\bm{u}(x_i,y_{j-\frac{1}{2}},t))}dt}{\Delta y}\\
&=\frac{1}{\Delta t}\int_{t_{n+1}}^{t_n}\frac{1}{4}[\bm{S}(x_{i+\frac{1}{2}},y_j,\bm{u}(x_{i+\frac{1}{2}},y_j,t))+\bm{S}(x_{i-\frac{1}{2}},y_j,\bm{u}(x_{i-\frac{1}{2}},y_j,t))+\\ 
&\quad\quad\quad\quad\quad\quad\;\bm{S}(x_i,y_{j+\frac{1}{2}},\bm{u}(x_i,y_{j+\frac{1}{2}},t))+\bm{S}(x_i,y_{j-\frac{1}{2}},\bm{u}(x_i,y_{j-\frac{1}{2}},t))]dt,\\
&n=0,1,\cdot\cdot\cdot,N-1;
\;i=1,2,\cdot\cdot\cdot,\frac{b-a}{\Delta x};
\;j=1,2,\cdot\cdot\cdot,\frac{d-c}{\Delta y},
\end{aligned}
\label{FiniteVolumeMethod}
\end{equation}
where 
\begin{equation}
\overline{\bm{u}}_{i,j}^n
=\frac{1}{\Delta x\Delta y}\int_{x_{i-\frac{1}{2}}}^{x_{i+\frac{1}{2}}}dx\int_{y_{j-\frac{1}{2}}}^{y_{j+\frac{1}{2}}}
\bm{u}(x,y,t_{n})dy,\ \
n=0,1,\cdot\cdot\cdot,N.
\end{equation}

If we approximate the flux with Riemann solver, the numerical error is measured by $\triangle\bm{u}$, i.e., the
jump across the interface, 
\begin{equation}
\text{numerical flux at} \; (x_{i+\frac{1}{2}},y_j)-\frac{1}{\Delta t_n}\int_{t_n}^{t_{n+1}}{\bm F}({\bm u}(x_{i+\frac{1}{2}},y_j,t))dt=\mathcal{O}(\lvert \triangle \bm{u}\lvert).
\label{NumericalErrorRiemann}
\end{equation}
When strong waves are present, this error in \eqref{NumericalErrorRiemann} will not decrease with the refinement of the grid, and will cause significant damage to the numerical simulation.
While if we approximate the flux with GRP solver, the numerical error is measured by $\Delta t$, which is equivalent to grid size,
\begin{equation}
\begin{aligned}
&\text{numerical flux at} \; (x_{i+\frac{1}{2}},y_{j})-\frac{1}{\Delta t_n}\int_{t_n}^{t_{n+1}}{\bm F}({\bm u}(x_{i+\frac{1}{2}},y_j,t))dt\\
&=-\frac{(\Delta t_n)^2}{6}\frac{\partial^2\bm{F}}{\partial t^2}(\bm{u}_{i+\frac{1}{2},j}^n)
+\mathcal{O}((\triangle t_{n})^3).
\end{aligned}
\label{NumericalErrorGRP}
\end{equation}
The accuracy of \eqref{NumericalErrorGRP} is $\mathcal{O}((\triangle t_{n})^3)$ for smooth solution,
and the accuracy of \eqref{NumericalErrorGRP} is $\mathcal{O}((\triangle t_{n})^2)$ for discontinuous solution. 
The temporal integral in \eqref{FiniteVolumeMethod} is with first-order and second-order accuracy respectively,
when we use the Riemann solver and the GRP solver. The other advantages of the GRP solver over the Riemann solver are:
1) the thermodynamic entropy changes are embedded; 
2) the effects of source term and tangential fluxes are embedded.  
We will elaborate on the above later.

\section{The approximation of the numerical fluxes}

Now we consider the approximation of the numerical flux at $(x_{i+\frac{1}{2}},y_j)$
for
\begin{equation} 
\frac{1}{\Delta t_n}\int_{t_n}^{t_{n+1}}{\bm F}({\bm u}(x_{i+\frac{1}{2}},y_j,t))dt
\label{Fluxes}
\end{equation}
with Riemann solver and GRP solver, respectively.

\subsection{The numerical fluxes with Riemann solver}

The numerical fluxes with Riemann solver consider the following Riemann problem along the normal direction of the cell interface
for Riemann solution, i.e., ${\bm u}_{i+\frac{1}{2},j}^{n,RP}$,
\begin{equation}
\begin{aligned}
&\partial_t \bm{u}(x,y_j,t) + \partial_x \bm{F}({\bm u}(x,y_j,t))=\bm{0},\ \ t>t_n,\\
&\bm{u}(x,y_j,t_n)=
\begin{cases}
\bm{u}_{i+1,j}^n,&x<x_{i+\frac{1}{2}},\\
\bm{u}_{i,j}^n,&x>x_{i+\frac{1}{2}}.
\end{cases}
\end{aligned}
\label{RP_RP}
\end{equation}
Then the numerical flux at $(x_{i+\frac{1}{2}},y_j)$ for \eqref{Fluxes} can be approximated as
\begin{equation}
\frac{1}{\Delta t_n}\int_{t_n}^{t_{n+1}}{\bm F}({\bm u}(x_{i+\frac{1}{2}},y_j,t))dt
\approx
{\bm F}\left({\bm u}_{i+\frac{1}{2},j}^{n,RP}\right).
\label{Flux_RP}
\end{equation}

\subsection{The numerical fluxes with GRP solver}
The numerical fluxes with GRP solver consider the following GRP along the normal direction of the cell interface for Generalized Riemann solution, i.e.,
${\bm u}_{i+\frac{1}{2},j}^{n+\frac{1}{2}}={\bm u}_{i+\frac{1}{2},j}^{n,GRP}+\frac{\Delta t_n}{2}\left(\frac{\partial \bm u}{\partial t}\right)_{i+\frac{1}{2},j}^{n}$,
\begin{equation}
\begin{aligned}
&\partial_t \bm{u}(x,y_j,t) + \partial_x \bm{F}({\bm u}(x,y_j,t))=- \partial_y \bm{G}({\bm u}(x,y_j,t))+\bm{S}(x,y_j,{\bm u}(x,y_j,t)),\ \ t>t_n,\\
&\bm{u}(x,y_j,t_n)=
\begin{cases}
\bm{u}_{i,j}^n+(x-x_i){\bm\sigma}_{i,j,x}^n,&x<x_{i+\frac{1}{2}},\\
\bm{u}_{i+1,j}^n+(x-x_{i+1}){\bm\sigma}_{i+1,j,x}^n,&x>x_{i+\frac{1}{2}}.
\end{cases}
\end{aligned}
\end{equation}
Since we just want to construct the flux normal to cell interfaces, the tangential effect can be regarded as part of the source.

Then the numerical flux at $(x_{i+\frac{1}{2}},y_j)$ for \eqref{Fluxes} can be approximated as
\begin{equation}
\frac{1}{\Delta t_n}\int_{t_n}^{t_{n+1}}{\bm F}({\bm u}(x_{i+\frac{1}{2}},y_j,t))dt
\approx
{\bm F}\left({\bm u}_{i+\frac{1}{2},j}^{n+\frac{1}{2}}\right).
\label{Flux_GRP}
\end{equation}

\subsubsection{Quasi 1-D acoustic case}
If $\lVert \bm{u}(x_{i+\frac{1}{2}}+0,y_j,t_n)-\bm{u}(x_{i+\frac{1}{2}}-0,y_j,t_n)\lVert$ is very small, 
which is  proportional to the mesh size, then we  view it as the quasi 1-D acoustic case. 

Denote the approximate Riemann solution by 
\begin{equation}
\bm{u}_{i+\frac{1}{2},j}^n:= \frac{1}{2}\left(\bm{u}(x_{i+\frac{1}{2}}+0,y_j,t_n)+ \bm{u}(x_{i+\frac{1}{2}}-0,y_j,t_n)\right),
\end{equation}
and 
\begin{equation}
A\left(\bm{u}_{i+\frac{1}{2},j}^n\right):=\frac{\partial \bm{F}}{\partial \bm{u}}\left(\bm{u}_{i+\frac{1}{2},j}^n\right).
\end{equation}
We make the decomposition 
\begin{equation}
A\left(\bm{u}_{i+\frac{1}{2},j}^n\right)=R\Lambda R^{-1},
\end{equation}
where $\Lambda =\text{diag}\{\lambda_{i}\}$, $R$ is the left eigenmatrix of $A\left(\bm{u}_{i+\frac{1}{2},j}^n\right)$.
Then the acoustic GRP solver takes
\begin{equation}
\begin{aligned}
&\left(\frac{\partial \bm{u}}{\partial t}\right)_{i+\frac{1}{2},j}^n
=-R\Lambda^{+}R^{-1}\left(\frac{\partial\bm{u}}{\partial x}\right)\left(x_{i+\frac{1}{2}}-0,y_j,t_n\right)
      -R\Lambda^{-}R^{-1}\left(\frac{\partial\bm{u}}{\partial x}\right)\left(x_{i+\frac{1}{2}}+0,y_j,t_n\right)\\
   &\quad\quad\;-RI^{+}R^{-1}\left(\frac{\partial \bm{G}}{\partial y}-S\right)\left(x_{i+\frac{1}{2}}-0,y_j,t_n\right)
      -RI^{-}R^{-1}\left(\frac{\partial \bm{G}}{\partial y}-S\right)\left(x_{i+\frac{1}{2}}+0,y_j,t_n\right),
\end{aligned}
\end{equation}
where $\Lambda^{+}=\text{diag}\{\max(\lambda_i,0)\}$,
$\Lambda^{-}=\text{diag}\{\min(\lambda_i,0)\}$,
$I^{+}=\frac{1}{2}\text{diag}\{1+\text{sign}(\lambda_i)\}$,
$I^{-}=\frac{1}{2}\text{diag}\{1-\text{sign}(\lambda_i)\}$.

\subsubsection{Quasi 1-D genuinely nonlinear case}
If $\lVert \bm{u}(x_{i+\frac{1}{2}}+0,y_j,t_n)-\bm{u}(x_{i+\frac{1}{2}}-0,y_j,t_n)\lVert$ is large, 
which is not proportional to the mesh size, then we view it as the quasi 1-D genuinely nonlinear case.

We solve the 1-D Riemann problem for 
\begin{equation}
\begin{aligned}
&\frac{\partial\bm{u}}{\partial t}+\frac{\partial \bm{F}(\bm{u})}{\partial x}={\bm 0},\quad t>0,\\
&\bm{u}(x,y_j,t_n)=
\begin{cases}
\bm{u}_{{i+\frac{1}{2}}^{-},j}^n,\\
\bm{u}_{{i+\frac{1}{2}}^{+},j}^n,
\end{cases}
\end{aligned}
\label{RP_GRP}
\end{equation} 
to obtain the local Riemann solution $\bm{u}_{i+\frac{1}{2},j}^n$,
where
\begin{equation}
\begin{aligned}
\bm{u}_{{i+\frac{1}{2}}^{-},j}^n:=\bm{u}_{i,j}^n+\frac{\Delta x}{2}{\bm\sigma}_{i,j,x}^n,\quad x<x_{i+\frac{1}{2}},\\
\bm{u}_{{i+\frac{1}{2}}^{+},j}^n:=\bm{u}_{i+1,j}^n-\frac{\Delta x}{2}{\bm\sigma}_{i+1,j,x}^n,\quad x>x_{i+\frac{1}{2}},
\end{aligned}
\label{RP_GRP_initialdata}
\end{equation}
Then we solve the 1-D GRP for 
\begin{equation}
\begin{aligned}
&\frac{\partial\bm{u}}{\partial t}+\frac{\partial \bm{F}(\bm{u})}{\partial x}=0,\quad t>0,\\
&\bm{u}(x,y_j,t_n)=
\begin{cases}
\bm{u}_{i+\frac{1}{2}^{-},j}^n+(x-x_{i+\frac{1}{2}}){\bm\sigma}_{i,j,x}^n,&x<x_{i+\frac{1}{2}},\\
\bm{u}_{i+\frac{1}{2}^{+},j}^n+(x-x_{i+\frac{1}{2}}){\bm\sigma}_{i+1,j,x}^n,&x>x_{i+\frac{1}{2}},
\end{cases}
\end{aligned}
\end{equation}
to obtain the local GRP solution $\left(\frac{\partial \bm{u}}{\partial t}\right)_{i+\frac{1}{2},j}^{n,x}$ by the 1-D GRP solver as \cite{Ben-Artzi-2006}.
Espectially for strong sparse wave, the GRP solver have
\begin{equation}
\frac{D u_x}{Dt}+\frac{Dp}{Dt}=\frac{D\psi}{Dt},
\label{entropyVariants}
\end{equation}
where $\frac{D}{Dt}=\frac{\partial}{\partial t}+{u_x}\frac{\partial}{\partial x}$,
$\frac{D\psi}{Dt}=cK(\rho,S)\frac{\partial S}{\partial x}-c\frac{\partial \psi}{\partial x}$,
$S$ is the entropy,
and $\psi$ is a Riemann variant.

We decompose $A(\bm{u}_{i+\frac{1}{2},j}^n)=\frac{\partial \bm{F}}{\partial \bm{u}}(\bm{u}_{i+\frac{1}{2},j}^n)=R\Lambda R^{-1}$,
$\Lambda =\text{diag}\{\lambda_{i}\}$, $R$ is the left eigenmatrix of $A\left(\bm{u}_{i+\frac{1}{2},j}^n\right)$.
Then we set
\begin{equation}
\bm{H}(x,y_j,t_n)=
\begin{cases}
-RI^{+}R^{-1}\left(\frac{\partial \bm{G}\left(\bm{u}_{{i+\frac{1}{2}}^{-},j}^n\right)}{\partial y}-{\bm S}(\bm{u}_{{i+\frac{1}{2}}^{-},j}^n)\right)_{(x_{i+\frac{1}{2}}-0,y_j)}
,& x<x_{i+\frac{1}{2}},\\
-RI^{-}R^{-1}\left(\frac{\partial \bm{G}\left(\bm{u}_{{i+\frac{1}{2}}^{+},j}^n\right)}{\partial y}-{\bm S}(\bm{u}_{{i+\frac{1}{2}}^{+},j}^n)\right)_{(x_{i+\frac{1}{2}}+0,y_j)}
, & x>x_{i+\frac{1}{2}},
\end{cases}
\label{SourceTermTangentialTerm}
\end{equation}
where $I^{+}=\frac{1}{2}\text{diag}\{1+\text{sign}(\lambda_i)\}$, $I^{-}=\frac{1}{2}\text{diag}\{1-\text{sign}(\lambda_i)\}$.
 
Then we solve the quasi 1-D GRP 
\begin{equation}
\frac{\partial\bm u}{\partial t}+\frac{\partial \bm{F}(\bm{u})}{\partial x}={\bm H}(x,y_j,t_n)
\end{equation}
to obtain $\left(\frac{\partial \bm{u}}{\partial t}\right)_{i+\frac{1}{2},j}^n$ approximately as
\begin{equation}
\left(\frac{\partial \bm{u}}{\partial t}\right)_{i+\frac{1}{2},j}^n
:=\left(\frac{\partial \bm{u}}{\partial t}\right)_{i+\frac{1}{2},j}^{n,x}+
{\bm H}(x,y_j,t_n).
\label{AddSourceTermTangentialTerm}
\end{equation}

Compared with the numerical fluxes with Riemann solver, the advantages of the numerical fluxes with GRP solver are as follows:
1) Having higher numerical accuracy, which can be reflected from \eqref{Flux_RP},\eqref{Flux_GRP},\eqref{RP_RP},\eqref{RP_GRP},\eqref{RP_GRP_initialdata};
2) Containing the effects of source term and tangential term, which can be reflected from \eqref{SourceTermTangentialTerm}, \eqref{AddSourceTermTangentialTerm};
3) Contains thermodynamic entropy change, which can be reflected from \eqref{entropyVariants}.

\section{The definition of the ghost fluid states}
In this section, we give the definition of the ghost fluid states, i.e., $\bm{u}_k^{g,n}$ for the Riemann problem based ghost fluid method
and the GRP based ghost fluid method. Then compare the performance differences between them.

When we transform the control equation of \eqref{MathModel} from the coordinate system $(x,y)$ to the coordinate system $(\xi,\eta)$, according to the invariance principle of the fundamental equations of fluid mechanics, we have
\begin{equation}
\frac{\partial \widetilde{\bm u}}{\partial t}+\frac{\partial \bm F(\widetilde{\bm u})}{\partial \xi}+\frac{\partial \bm G(\widetilde{\bm u})}{\partial \eta}={\bm S}(\xi,\eta,{\widetilde{\bm u}}).
\end{equation}

\subsection{The Riemann problem based ghost fluid method}
The ghost fluid states, i.e., $\bm{u}_k^{g,n}$ for the Riemann problem based ghost fluid method is defined by
establishing the following double-medium Riemann problem,
 \begin{equation}
\begin{aligned}
&\frac{\partial \widetilde{\bm u}(\xi,0,t)}{\partial t}+\frac{\partial \bm F(\widetilde{\bm u}(\xi,0,t))}{\partial \xi}={\bm 0}, \quad t>t^n,\\
&{\widetilde{\bm u}}(\xi,0,t^n)=
\begin{cases}
\widetilde{\bm u}_{1},&\xi<0,\\
\widetilde{\bm u}_{2},&\xi>0,
\end{cases}
\end{aligned}
\label{doublemediumRP}
\end{equation}
where 
\begin{equation}
{\widetilde{\bm u}}_{k}=(\rho_{k}, u_{\xi,k}, p_{k})^{\top}, \quad u_{\xi,k}=u_{x,k}n_x+u_{y,k}n_y, \quad k=1,2,
\end{equation}
and 
\begin{equation}
{\bm u}_{k}=(\rho_k,u_{x,k},u_{y,k},p_k)^{\top}, \quad k=1,2,
\end{equation}
 are obtained by the weighted least square fitting as \cite{Huo-2024}.

 Solve the problem  \eqref{doublemediumRP} for the after wave states of each medium, i.e., 
\begin{equation}
{\widetilde{\bm u}}^{*}_{k}=(\rho^{*}_{k},u_{\xi,k}^{*},p_{k}^{*})^{\top}, \quad k=1,2,
\end{equation}
 based on the fact that the velocities are equal and the pressures are balanced across the interface.
The densities are obtained by the after wave states seperately as \cite{Huo-2024}.

Then the Riemann problem based ghost fluid states are defined as
\begin{equation}
{\bm u}_k^{g,n}(x,y):={\bm u}^{*}_{k}=(\rho_{k}^{*},u_{\xi,k}^{*}n_x+u_{\eta,k}t_x,u_{\xi,k}^{*}n_y+u_{\eta,k}t_y,p_{k}^{*})^{\top}, \quad k=1,2,
\end{equation}
where the unit tangential vector of the multimedium interface is $T=(t_x,t_y)^{\top}=(-n_y,n_x)^{\top}$.

\subsection{The GRP based ghost fluid method}
The ghost fluid states, i.e., $\bm{u}_k^{g,n}$ for the GRP based ghost fluid method is defined by
establishing a double-medium GRP along the normal direction of the material interface.
Since the governing equations of computational fluid dynamics satisfy Galilean invariance \cite{Piekarski-2007,Robert-2010}, their specific form is given as follows:
\begin{equation}
\begin{aligned}
&\frac{\partial {\widetilde{\bm u}}(\xi,0,t)}{\partial t}+\frac{\partial \bm F({\widetilde{\bm u}}(\xi,0,t))}{\partial \xi}
=-\frac{\partial \bm G({\widetilde{\bm u}}(\xi,0,t))}{\partial \eta}+{\widetilde{\bm S}}(\xi,0,{{\widetilde{\bm u}}(\xi,0,t)}),\quad t>t^n,\\
&{\widetilde{\bm u}}(\xi,0,t^n)=
\begin{cases}
{\widetilde{\bm u}}_{1}+\xi\left(\frac{\partial {\widetilde{\bm u}}}{\partial \xi}\right)_1,&\xi<0,\\
{\widetilde{\bm u}}_{2}+\xi\left(\frac{\partial {\widetilde{\bm u}}}{\partial \xi}\right)_2,&\xi>0,
\end{cases}
\end{aligned}
\label{doublemediumGRP}
\end{equation}
where $\xi$ and $\eta$ are the normal and tangential directions of the material interface,
\begin{equation}
\left(\frac{\partial {\widetilde{\bm u}}}{\partial \xi}\right)_k
=
\left[
\begin{array}{c}
\frac{\partial \rho}{\partial \xi}\\[3pt] 
\frac{\partial u_{\xi}}{\partial \xi}\\[3pt]
\frac{\partial u_{\eta}}{\partial \xi}\\[3pt]
\frac{\partial p}{\partial \xi}
\end{array}
\right]_k
=
\left[
\begin{array}{c}
\frac{\partial \rho}{\partial x}\\[3pt]
\frac{\partial u_x}{\partial x}n_x+\frac{\partial u_y}{\partial x}n_y\\[3pt]
\frac{\partial u_x}{\partial x}t_x+\frac{\partial u_y}{\partial x}t_y\\[3pt]
\frac{\partial p}{\partial x}
\end{array}
\right]_kn_x
+
\left[
\begin{array}{c}
\frac{\partial \rho}{\partial y}\\[3pt]
\frac{\partial u_x}{\partial y}n_x+\frac{\partial u_y}{\partial y}n_y\\[3pt]
\frac{\partial u_x}{\partial y}t_x+\frac{\partial u_y}{\partial y}t_y\\[3pt]
\frac{\partial p}{\partial y}
\end{array}
\right]_kn_y,\quad k=1,2,
\end{equation}
and 
\begin{equation}
\begin{aligned}
&\left(\frac{\partial \bm u}{\partial x}\right)_{k}=
\left(
\left(\frac{\partial \rho}{\partial x}\right)_k,
\left(\frac{\partial u_x}{\partial x}\right)_k,
\left(\frac{\partial u_y}{\partial x}\right)_k,
\left(\frac{\partial p}{\partial x}\right)_k \right)^{\top},\\
&\left(\frac{\partial \bm u}{\partial y}\right)_{k}=
\left(
\left(\frac{\partial \rho}{\partial y}\right)_k,
\left(\frac{\partial u_x}{\partial y}\right)_k,
\left(\frac{\partial u_y}{\partial y}\right)_k,
\left(\frac{\partial p}{\partial y}\right)_k\right)^{\top},\\
&k=1,2,
\end{aligned}
\end{equation}
are obtained by the weighted least square fitting as \cite{Huo-2024}. Solve the problem  \eqref{doublemediumGRP} for the after wave states of each medium, i.e., 
\begin{equation}
\left(\frac{\partial {\widetilde{\bm u}}}{\partial t}\right)_k^{*}=
\left(
\left(\frac{\partial \rho}{\partial t}\right)_k^*,
\left(\frac{\partial u_{\xi}}{\partial t}\right)_k^*,
\left(\frac{\partial u_{\eta}}{\partial t}\right)_k^*,
\left(\frac{\partial p}{\partial t}\right)_k^*
\right)^{\top},
\quad k=1,2,
\label{TemperalDerivatives}
\end{equation}
based on the fact that  
the velocities are equal and the pressures are balanced across the material interface
, then the material derivatives of velocity and pressure  are continuous across the material interface.
The temperal derivatives of densities are obtained by the after wave states of the double-medium GRP seperately.  

By combing the temperal derivatives \eqref{TemperalDerivatives} and the equations\eqref{doublemediumGRP},
we obtain the spatial derivatives of the after wave states for each medium, i.e.,
\begin{equation}
\left(\frac{\partial {\widetilde{\bm u}}}{\partial \xi}\right)_k^{*}=
\left(
\left(\frac{\partial \rho}{\partial \xi}\right)_k^*,
\left(\frac{\partial u_{\xi}}{\partial \xi}\right)_k^*,
\left(\frac{\partial u_{\eta}}{\partial \xi}\right)_k^*,
\left(\frac{\partial p}{\partial \xi}\right)_k^*
\right)^{\top},
\end{equation}
and
\begin{equation}
\begin{aligned}
\left(\frac{\partial {\widetilde{\bm u}}}{\partial \eta}\right)_k^{*}
&=
\left(
\left(\frac{\partial \rho}{\partial \eta}\right)_k^*,
\left(\frac{\partial u_{\xi}}{\partial \eta}\right)_k^*,
\left(\frac{\partial u_{\eta}}{\partial \eta}\right)_k^*,
\left(\frac{\partial p}{\partial \eta}\right)_k^*
\right)^{\top}\\
&=
\left(
\left(\frac{\partial \rho}{\partial \eta}\right)_k,
\left(\frac{\partial u_{\xi}}{\partial \eta}\right)_k,
\left(\frac{\partial u_{\eta}}{\partial \eta}\right)_k,
\left(\frac{\partial p}{\partial \eta}\right)_k
\right)^{\top}\\
&=\left[
\begin{array}{c}
\frac{\partial \rho}{\partial x}\\[3pt]
\frac{\partial u_x}{\partial x}n_x+\frac{\partial u_y}{\partial x}n_y\\[3pt]
\frac{\partial u_x}{\partial x}t_x+\frac{\partial u_y}{\partial x}t_y\\[3pt]
\frac{\partial p}{\partial x}
\end{array}
\right]_kt_x
+
\left[
\begin{array}{c}
\frac{\partial \rho}{\partial y}\\[3pt]
\frac{\partial u_x}{\partial y}n_x+\frac{\partial u_y}{\partial y}n_y\\[3pt]
\frac{\partial u_x}{\partial y}t_x+\frac{\partial u_y}{\partial y}t_y\\[3pt]
\frac{\partial p}{\partial y}
\end{array}
\right]_kt_y
\\
\end{aligned}
\end{equation}
\begin{equation}
k=1,2.
\end{equation}

Then the GRP based ghost fluid states are  defined as
\begin{equation}
{\bm u}_{k}^{g,n}(x,y):={\bm u}_k^{*}+(x-x_{\Gamma})\left(\frac{\partial {\bm u}}{\partial x}\right)_{k}^{*}+(y-y_{\Gamma})\left(\frac{\partial {\bm u}}{\partial y}\right)_{k}^{*},\quad k=1,2,
\end{equation}
where
\begin{equation}
\begin{aligned}
&\left(\frac{\partial {\bm u}}{\partial x}\right)_{k}^*=
\left[
\begin{array}{c}
\frac{\partial \rho}{\partial x}\\[3pt]
\frac{\partial u_{x}}{\partial x}\\[3pt]
\frac{\partial u_{y}}{\partial x}\\[3pt]
\frac{\partial p}{\partial x}
\end{array}
\right]_k^*
=
\left[
\begin{array}{c}
\frac{\partial \rho}{\partial \xi}\\[3pt]
\frac{\partial u_{\xi}}{\partial \xi}n_x+\frac{\partial u_{\eta}}{\partial \xi}t_x\\[3pt]
\frac{\partial u_{\xi}}{\partial \xi}n_y+\frac{\partial u_{\eta}}{\partial \xi}t_y\\[3pt]
\frac{\partial p}{\partial \xi}
\end{array}
\right]_k^*n_x
+
\left[
\begin{array}{c}
\frac{\partial \rho}{\partial \eta}\\[3pt]
\frac{\partial u_{\xi}}{\partial \eta}n_x+\frac{\partial u_{\eta}}{\partial \eta}t_x\\[3pt]
\frac{\partial u_{\xi}}{\partial \eta}n_y+\frac{\partial u_{\eta}}{\partial \eta}t_y\\[3pt]
\frac{\partial p}{\partial \eta}
\end{array}
\right]_k^*t_x,\\
&\left(\frac{\partial {\bm u}}{\partial y}\right)_{k}^*
=
\left[
\begin{array}{c}
\frac{\partial \rho}{\partial x}\\[3pt]
\frac{\partial u_{x}}{\partial x}\\[3pt]
\frac{\partial u_{y}}{\partial x}\\[3pt]
\frac{\partial p}{\partial x}
\end{array}
\right]_k^*
=
\left[
\begin{array}{c}
\frac{\partial \rho}{\partial \xi}\\[3pt]
\frac{\partial u_{\xi}}{\partial \xi}n_x+\frac{\partial u_{\eta}}{\partial \xi}t_x\\[3pt]
\frac{\partial u_{\xi}}{\partial \xi}n_y+\frac{\partial u_{\eta}}{\partial \xi}t_y\\[3pt]
\frac{\partial p}{\partial \xi}
\end{array}
\right]_k^*n_y
+
\left[
\begin{array}{c}
\frac{\partial \rho}{\partial \eta}\\[3pt]
\frac{\partial u_{\xi}}{\partial \eta}n_x+\frac{\partial u_{\eta}}{\partial \eta}t_x\\[3pt]
\frac{\partial u_{\xi}}{\partial \eta}n_y+\frac{\partial u_{\eta}}{\partial \eta}t_y\\[3pt]
\frac{\partial p}{\partial \eta}
\end{array}
\right]_k^*t_y,\\
&k=1,2.
\end{aligned}
\end{equation}

\section{Numerical tests}

\begin{figure}[H]
  \centering
  \includegraphics[width=12cm]{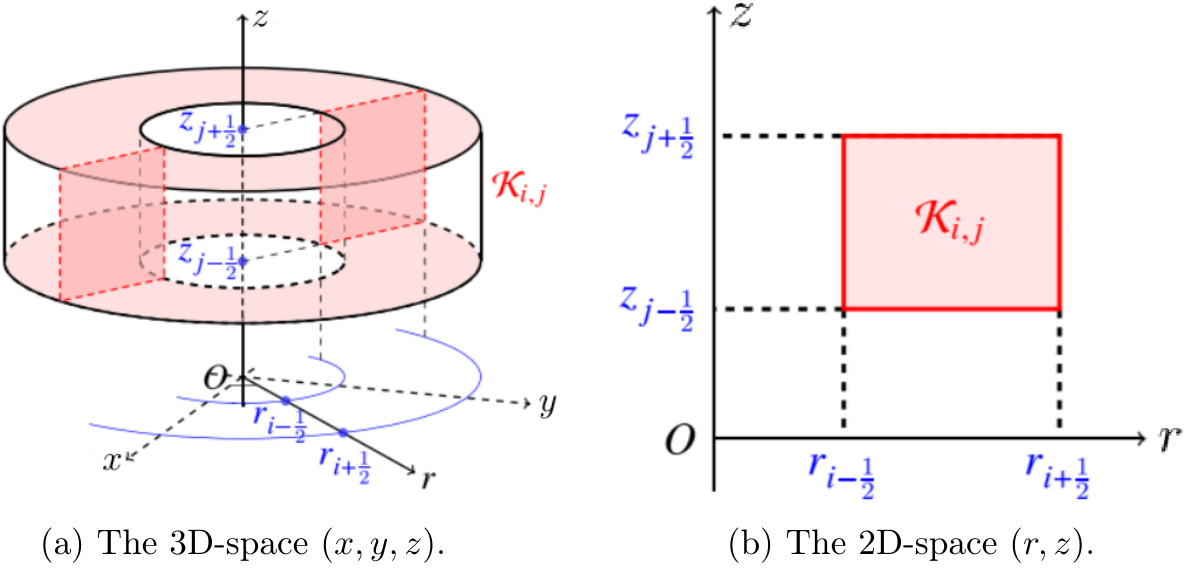}\\
  \caption{The schematic of an axisymmetric object and its cross section.}
  \label{Axisymmetric}
\end{figure}

Here we validate and compare the performance of the proposed algorithm for compressible multi-material flows with axisymmetric source terms in two-dimensional space. 
The cylindrical coordinates 
as shown in Figure \ref{Axisymmetric},
the Euler equations of axisymmtric fluid flows read
\begin{equation}
\frac{\partial \textbf{\emph{U}}}{\partial t}
+\frac{1}{r}\frac{\partial(r\textbf{\emph{F}}(\textbf{\emph{U}}))}{\partial r}
+\frac{\partial \textbf{\emph{G}}(\textbf{\emph{U}})}{\partial z}
=\frac{\textbf{\emph{H}}(\textbf{\emph{U}})}{r},
\label{AxisymmetricFlows}
\end{equation}
\begin{equation}
\textbf{\emph{H}}(\textbf{\emph{U}})
=\left[
0,
p,
0,
0
\right]^{\top}.
\end{equation}
The system \eqref{AxisymmetricFlows} is equivalent to 
\begin{equation}
\frac{\partial \textbf{\emph{U}}}{\partial t}+
\frac{\partial (\textbf{\emph{F}}(\textbf{\emph{U}}))}{\partial r}+
\frac{\partial \textbf{\emph{G}}(\textbf{\emph{U}})}{\partial z}=
\textbf{\emph{S}}(r,\textbf{\emph{U}}),
\end{equation}
\begin{equation}
\textbf{\emph{S}}(r,\textbf{\emph{U}})=
-\frac{1}{r}\left[\rho u, \rho u^2, \rho u v, u(\rho E+p)\right]^{\top}.
\end{equation}
Then the mathematical model of compressible multi-medium 
axisymmetric fluid flows is
\begin{equation}
\begin{aligned}
&\frac{\partial \textbf{\emph{U}}}{\partial t}+
\frac{\partial (\textbf{\emph{F}}(\textbf{\emph{U}}))}{\partial r}+
\frac{\partial \textbf{\emph{G}}(\textbf{\emph{U}})}{\partial z}=
\textbf{\emph{S}}(r,\textbf{\emph{U}}),\;\;\text{for}\;(r,z)\in\Omega,\; t>0,\\
&\textbf{\emph{U}}=
\begin{cases}
\textbf{\emph{U}}_1(r,z),\;\;\text{if}\;(r,z)\in\Omega_1(t),\\
\textbf{\emph{U}}_2(r,z),\;\;\text{if}\;(r,z)\in\Omega_2(t),
\end{cases}
e=
\begin{cases}
e_1(\rho,p),\;\;\text{if}\;(r,z)\in\Omega_1(t),\\
e_2(\rho,p),\;\;\text{if}\;(r,z)\in\Omega_2(t).
\end{cases}
\end{aligned}
\label{CMMAFF}
\end{equation}
The materials are modeled using the stiffened gas equation of state (EOS), expressed as:
\begin{equation}
e_{k}=\frac{p_k+\gamma_k p_{\infty}^{k}}{(\gamma_k-1)\rho_k}, \quad k=1,2.
\label{EOS_k}
\end{equation}

\subsection{Spherical bubble shock interaction problem}

\begin{figure}[H]
   \centering
   \includegraphics[width=13cm]{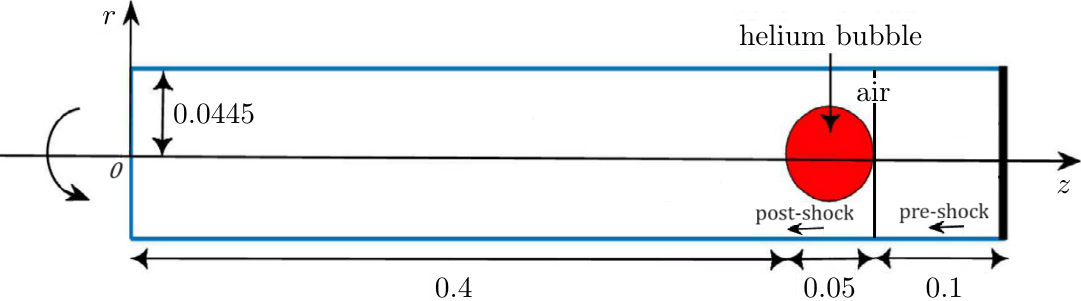}\\
   \caption{Schematic of initial flow configuration for 'spherical bubble shock interaction problem'.}
   \label{Test_initial}
\end{figure}

\begin{figure}[H]
   \centering
   \includegraphics[width=15cm]{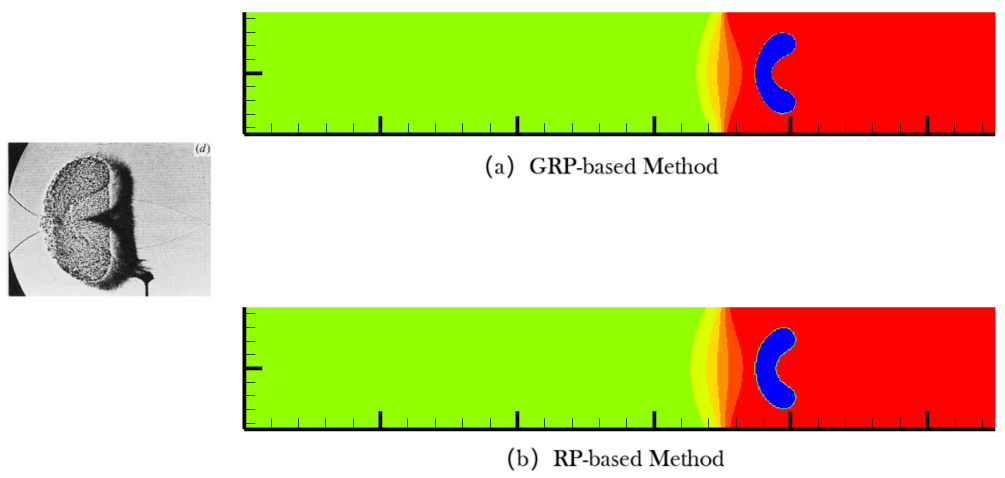}\\
   \caption{Shadow-photographs of the interaction of an $M_s=1.25$ shock wave moving from
right to left over a spherical helium volume (4.5 cm diameter) at $t_4=223\mu s$.}
   \label{Test_results_4}
\end{figure}

\begin{figure}[H]
   \centering
   \includegraphics[width=15cm]{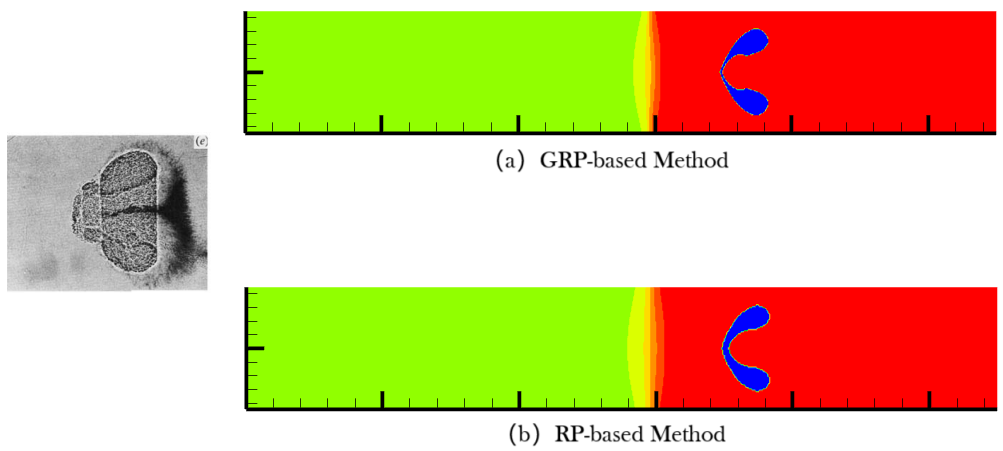}\\
   \caption{Shadow-photographs of the interaction of an $M_s=1.25$ shock wave moving from
right to left over a spherical helium volume (4.5 cm diameter) at $t_5=350\mu s$.}
   \label{Test_results_5}
\end{figure}

\begin{figure}[H]
   \centering
   \includegraphics[width=15cm]{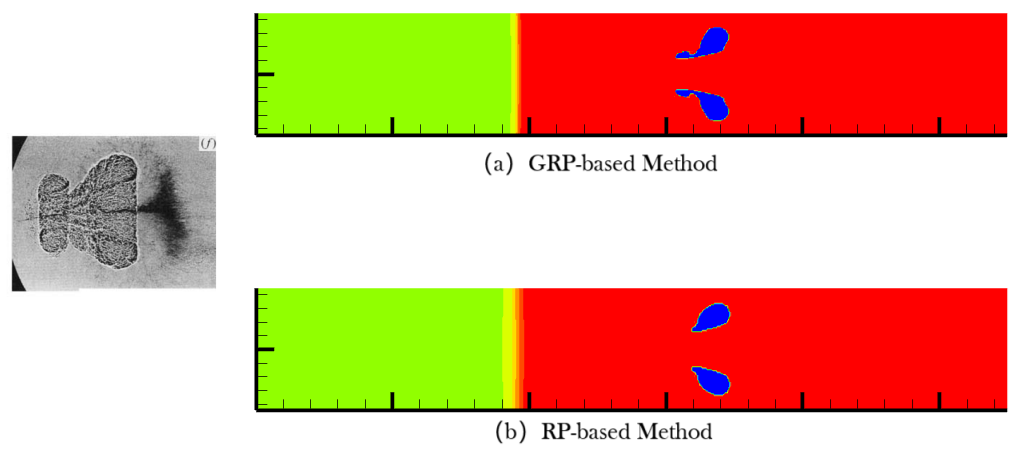}\\
   \caption{Shadow-photographs of the interaction of an $M_s=1.25$ shock wave moving from
right to left over a spherical helium volume (4.5 cm diameter) at $t_6=600\mu s$.}
   \label{Test_results_6}
\end{figure}

\begin{figure}[H]
   \centering
   \includegraphics[width=15cm]{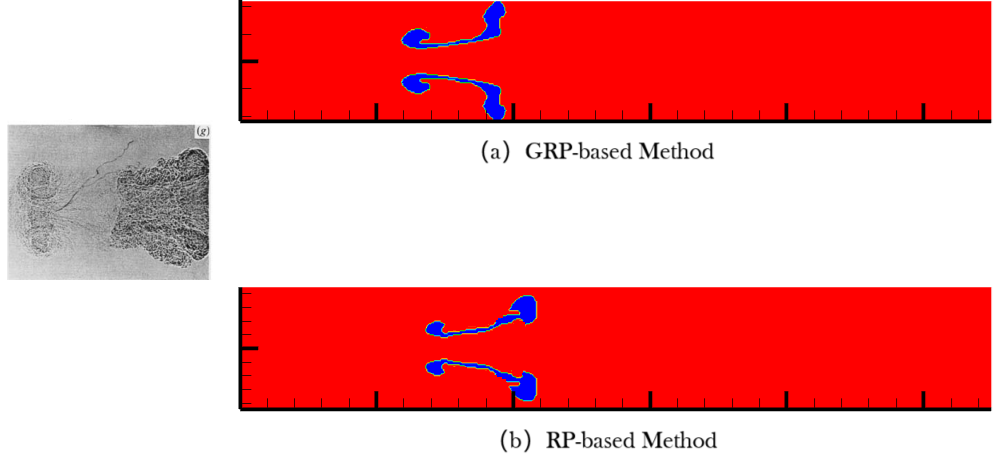}\\
   \caption{Shadow-photographs of the interaction of an $M_s=1.25$ shock wave moving from
right to left over a spherical helium volume (4.5 cm diameter) at $t_{\text{final}}=1594\mu s$.}
   \label{Test_results_end}
\end{figure}

This experiment is proposed in \cite{Hass-1987}, which describes that a spherical helium bubble surrounded by air is at rest for the initial moment, and a leftward incident shock wave triggers the motion. 
As shown in Figure \ref{Test_initial}, the computation domain is $[0,0.55]\times[-0.0445,0.0445]$, 
which is divided uniformly into $550\times90$ cells. 
The center of the bubble is located at $(0.425,0)$  and its radius is 0.025. The top, botton, left boundaries are given the solid wall condition, and the right boundary is given the piston-like boundary condition, moving leftward with the velocity 
$u^{*}$. The Mach number of the shock wave is $M_s=1.25$.
The helium and air are treated as an ideal gas with polytropic indices $\gamma_1=1.648$, $\gamma_2=1.4$ in \eqref{EOS_k}, respectively. The initial data for helium are
\begin{equation}
(\rho_1,u_{r1},u_{z1},p_1)=(0.2163,0,0,10^5),\quad \text{for}\; (r-0.425)^2+z^2\leqslant0.025^2. 
\end{equation} 
The initial data for air are
\begin{equation}
\begin{aligned}
&(\rho_2,u_{r2},u_{z2},p_2)\\
&=\begin{cases}
(1.189,0,0,10^5), &\text{for}\;z\leqslant0.45, (r-0.425)^2+z^2>0.025^2\;\text{(pre-shock)}, \\
(1.6985715,-128.67802,0,1.65625\times10^5),&\text{for}\;z>0.45, (r-0.425)^2+z^2>0.025^2\;\text{(post-shock)}.
\end{cases}
\end{aligned}
\end{equation}
Then according to the Rankine-Hugoniot relations, we can obtain the velocity of the piston $u^{*}=-128.678$.  
We compare our simulation results for GRP-based Method and the RP-based Method with the test results  in \cite{Hass-1987} at the intermediate time 
$t_4=223\mu s$, 
$t_5=350\mu s$, 
$t_6=600\mu s$,
and the final time 
$t_{\text{final}}=1594\mu s$ in Figure \ref{Test_results_4}, Figure \ref{Test_results_5}, Figure \ref{Test_results_6} and Figure \ref{Test_results_end}, respectively.

\subsection{The spherical shock induced collapse of an air bubble in water}
\begin{figure}[H]
\centering
\includegraphics[width=8cm]{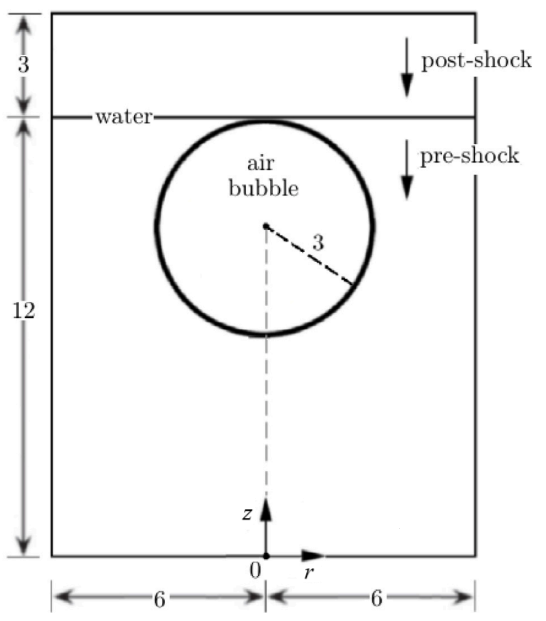}\\
\caption{Schematic of initial flow configuration for 'the  spherical  shock induced collapse of an air bubble in water'.}
\label{InitialAirBubble}
\end{figure}

\begin{figure}[H]
   \centering
   \includegraphics[width=16cm]{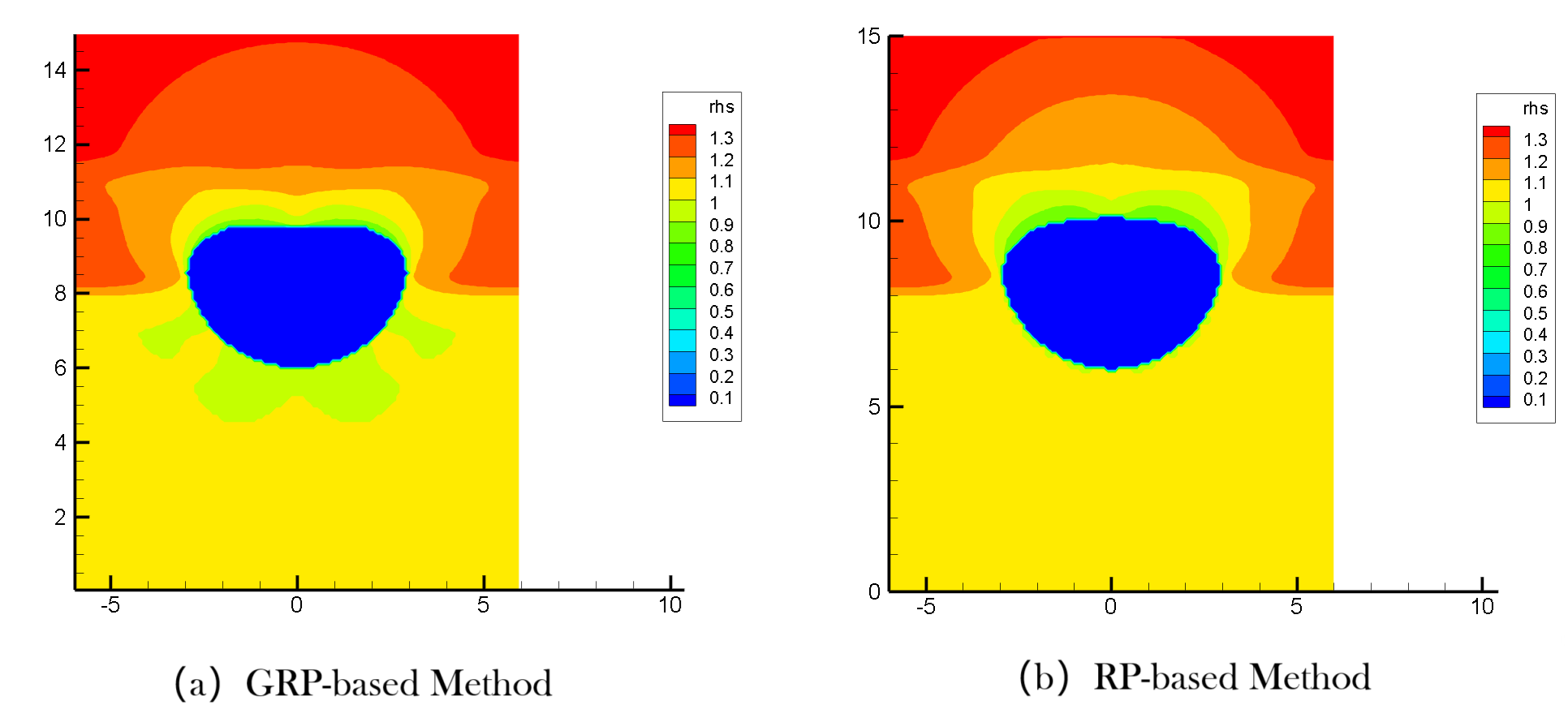}\\
   \caption{Simulation results of spherical underwater explosion 
at $t_1=0.012 \mu s$.
(a)GRP-based Method; (b)RP-based Method.
}
   \label{Test_results_t2_1}
\end{figure}

\begin{figure}[H]
   \centering
   \includegraphics[width=16cm]{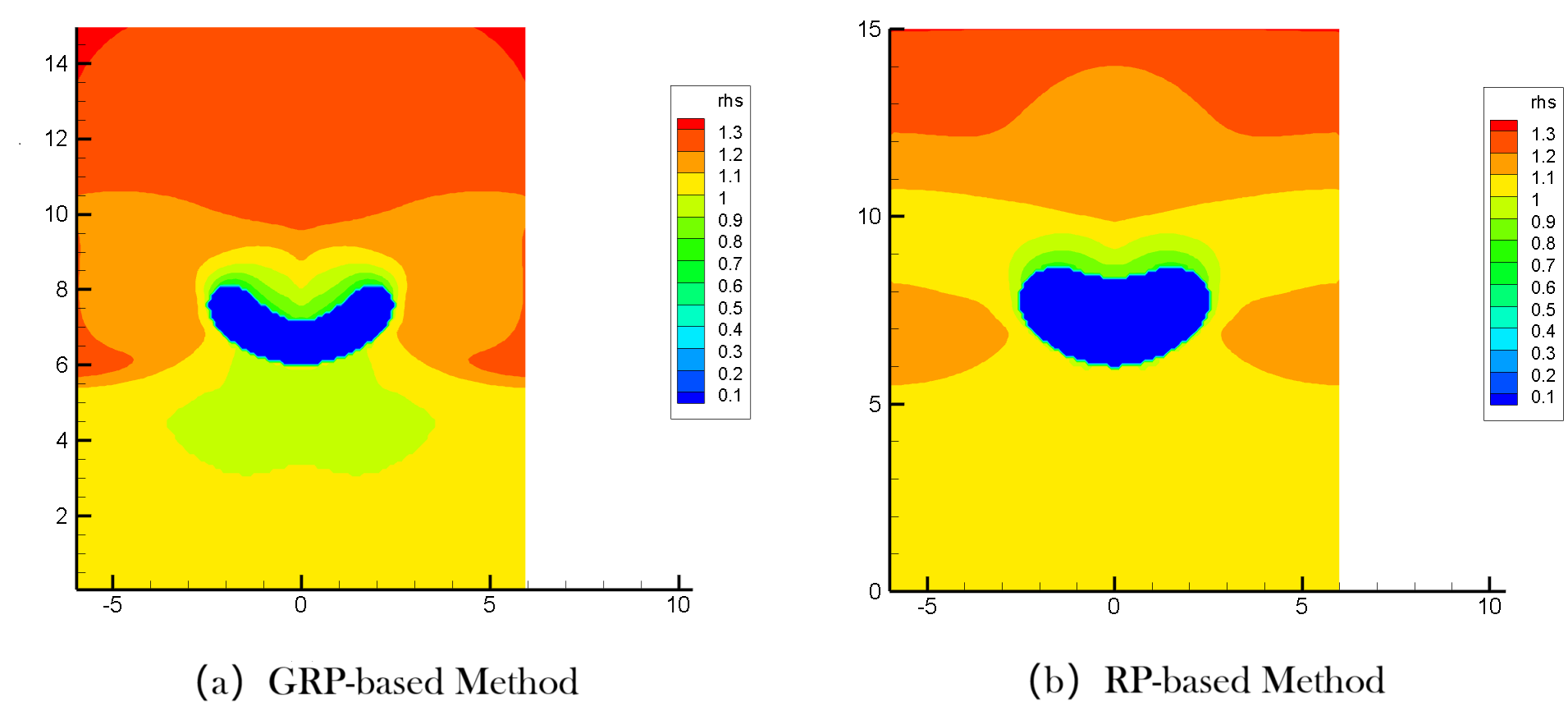}\\
   \caption{Simulation results of spherical underwater explosion 
at $t_2=0.020\mu s$.
(a)GRP-based Method; (b)RP-based Method.
}
   \label{Test_results_t2_2}
\end{figure}

\begin{figure}[H]
   \centering
   \includegraphics[width=16cm]{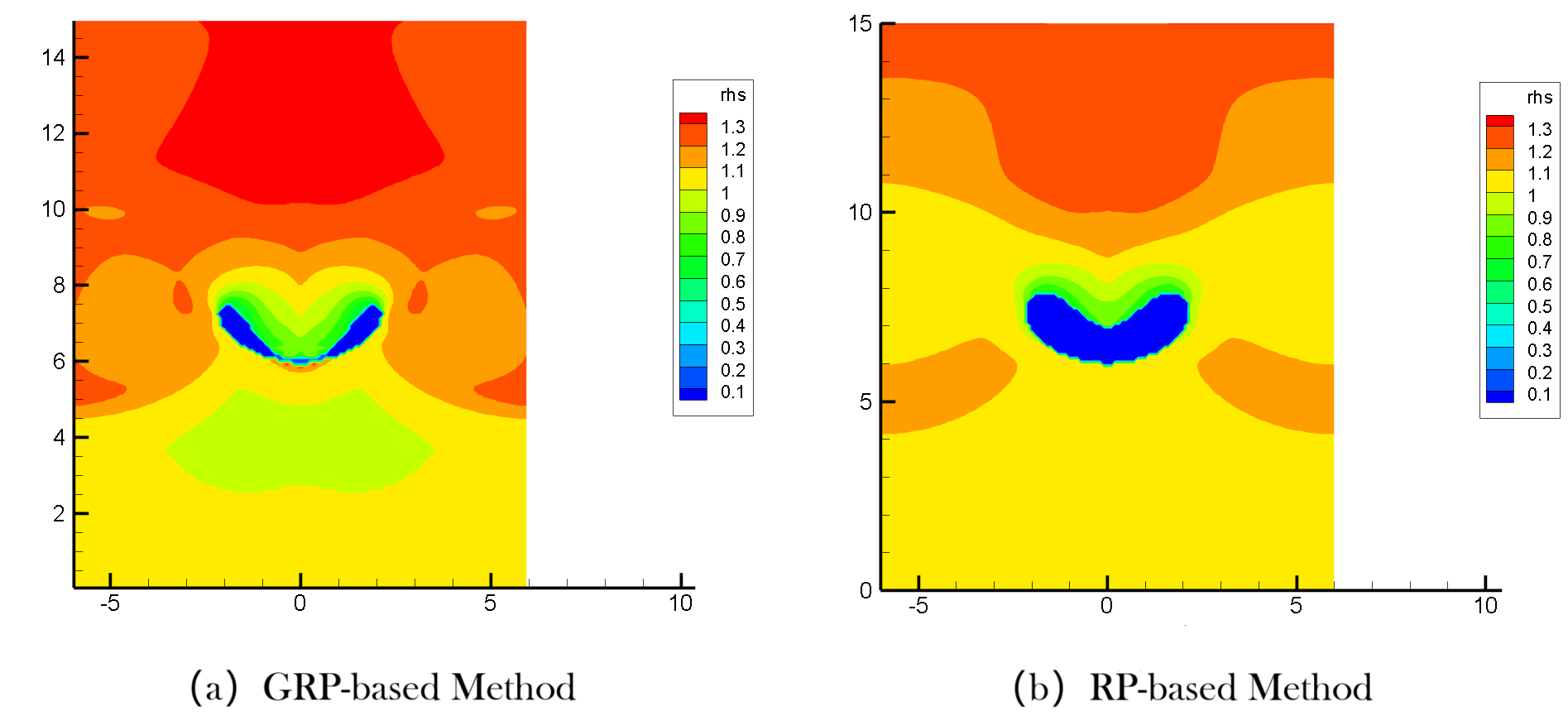}\\
   \caption{Simulation results of spherical underwater explosion 
at $t_3=0.02298\mu s$.
(a)GRP-based Method; (b)RP-based Method.
}
   \label{Test_results_t2_3}
\end{figure}

\begin{figure}[H]
   \centering
   \includegraphics[width=16cm]{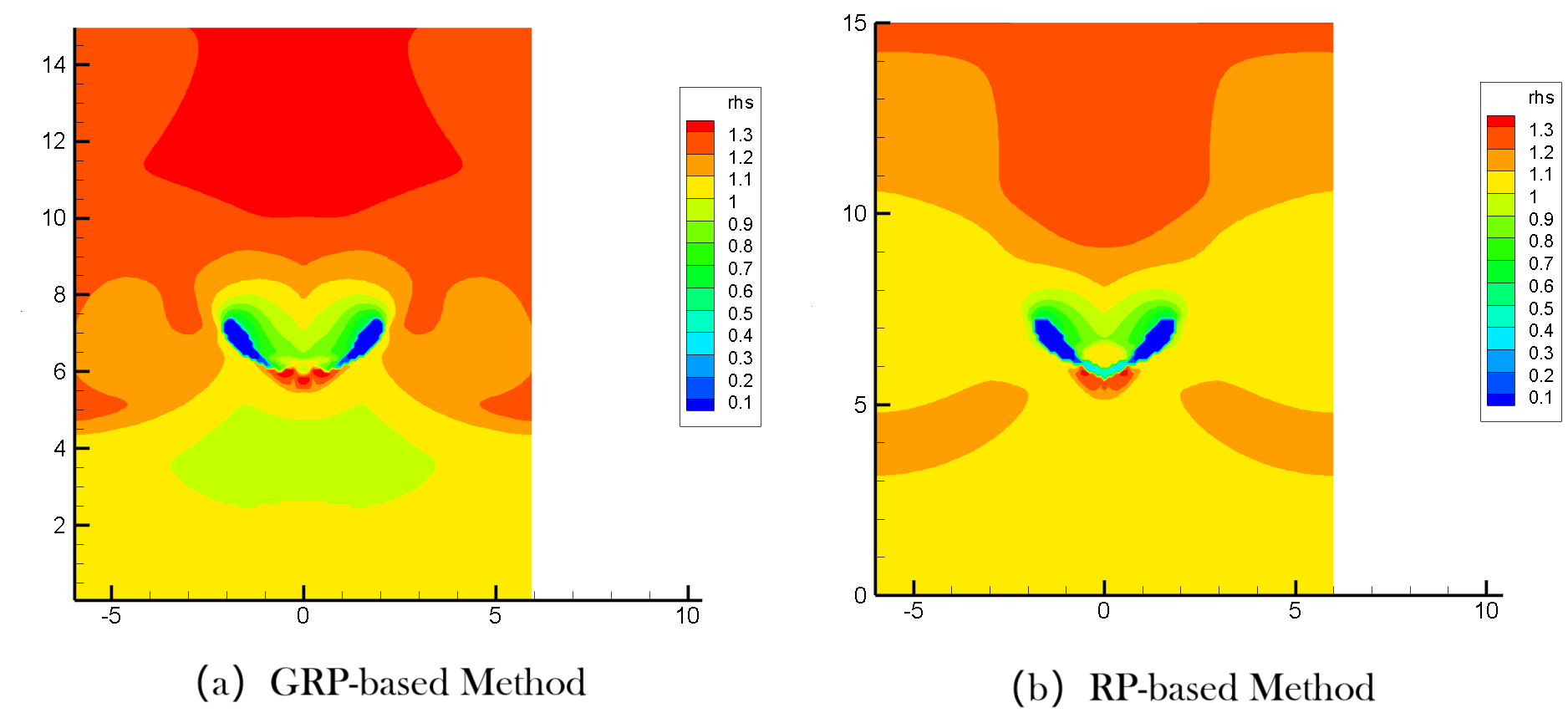}\\
   \caption{Simulation results of spherical underwater explosion 
at $t_4=0.02342\mu s$.
(a)GRP-based Method; (b)RP-based Method.
}
   \label{Test_results_t2_4}
\end{figure}

This problem is an experiment reported by Bourne and Field \cite{Bourne-1992},
and has been simulated by Ball et al. in \cite{Ball-2000}. It is important to some
practical applications. As illustrated in Fig. \ref{InitialAirBubble}, an speherical air bubble collapse
in water because of a shock wave propagates from up to down. The diameter of
the air bubble is 6. The compuational domain is $[-6,6]\times[0,15]$, which is divided
into $120\times150$ uniform cells. The initial states for air are
\begin{equation}
(\rho_1,u_{r1},u_{z1},p_1)
=(0.001,0,0,1),\quad\text{for}\; r^2+(z-9)^2\leqslant 9. 
\end{equation}
The shock is of pressure $19000$, then we can obtain the post-shock states by Rankine-Hugoniot relations. The initial states for water are
\begin{equation}
\begin{aligned}
&(\rho_2,u_{r2},u_{z2},p_2)\\
&=
\begin{cases}
(1,0,0,1), & \text{for}\; z\leqslant 12, r^2+(z-9)^2>9\;\text{(pre-shock)},\\
(1.313345,67.3267,0,19000), & \text{for}\; z>12, r^2+(z-9)^2>9\;\text{(post-shock)},\;
\end{cases}
\end{aligned}
\end{equation}
where $y=12$ is the localtion at which the shock wave hits the bubble. 
Figure \ref{Test_results_t2_1}, Figure \ref{Test_results_t2_2}, Figure \ref{Test_results_t2_3} and Figure \ref{Test_results_t2_4}
display the contour maps of  density of the simulation results for the GRP-based Method and RP-based Method at
the intermediate time $t=0.012\mu s$, $t=0.020 \mu s$, $t=0.02298\mu s$, and the terminal time $t=0.02342\mu s$, respectively. 

\section{Conclusions and future work}
The numerical results demonstrate that employing the generalized Riemann solver for numerical fluxes construction and ghost fluid states definition yields significantly more accurate and refined simulation outcomes compared to the conventional Riemann solver approach.
The fundamental reason lies in the fact that, compared to conventional Riemann solvers, the generalized Riemann solver accounts for the effects of source terms, tangential fluxes, and thermodynamic mechanisms. This leads to essential differences in problems involving source terms, multidimensional spaces, and compressible flows. In contrast, methods based on standard Riemann problems neglect these critical factors, resulting in loss of physical fidelity and diminished accuracy and reliability in numerical simulations.

Future work will focus on extending the proposed method to numerical simulations of multi-material reactive flows. A comprehensive comparison will be conducted to evaluate its advantages over conventional Riemann-based approaches, particularly in capturing reaction source terms, tangential effects, and thermodynamical mechanisms.
\section{Acknowledgments}

Zhixin Huo's research is supported by the Doctoral Fund of Henan Polytechnic University (Grant No. B2024-60),
and the 'double first-class' discipline creation project of surveying and mapping science and technology, Henan 
(Grant No. BSJJ202306).

\end{document}